# Optimal Shelter Location-Allocation during Evacuation with Uncertainties: A Scenario-Based Approach

Fardad Haghpanah[*] and Hamid Foroughi

*Department of Civil Engineering, Johns Hopkins University, Baltimore, MD 21218, USA*

[*] *Corresponding author, email: haghpanah@jhu.edu*

**Abstract**

Evacuation planning is an important and challenging element in emergency management due to the high level of uncertainty and numerous players and agencies involved in the event. To address all the factors with conflicting objectives, mathematical modeling has gained an extensive application over all aspects of evacuation planning to help responders and policy makers evaluate required time for evacuation and estimate numbers and distribution of casualties under different disaster scenarios. Correspondingly, mathematical formulation of evacuation optimization problems and solution methods are important when planning for evacuation. In this paper, the bi-level programming formulation of shelter location-allocation problem is considered. To account for stochasticity, a scenario-based approach is taken to address the uncertainty in the population to be evacuated from a small town in Lombardy region, Italy. Genetic algorithm is used as the solution method. Four scenarios are considered to study the optimal number and location of shelters for evacuation during normal weekdays, at nights, during weekends, and during vacation times with visiting travelers. The results highlight how different scenarios need different number and location of shelters for an optimal evacuation.

*Keywords: optimization, bi-level programming, genetic algorithm, evacuation*

## 1. Introduction

In many disasters, to ensure public safety and mitigate losses, it could be necessary to evacuate buildings, neighborhoods, or urban regions within the impact area. Accordingly, evacuation planning is an important element in emergency management. The main body of research on evacuation planning has been devoted to different subjects, such as evacuation behavior [1,2], traffic management [3,4], traffic emergency management policy [5,6], and origin-destination trip estimation [7,8]. What make emergency logistics more complex than regular business logistics are [9,10]: (1) additional uncertainties due to unusable roads, changing facilities' capacity, changing demand, and safety measures; (2) communication and coordination challenges due to damages, increasing number of involved parties, and lack of accurate real-time information; (3) critical need for timely and efficient delivery of services; and (4) limited resources, such as supply services, transportation and shelter capacity, fuel, and responders. All these factors bring conflicting objectives to the problem of optimal evacuation planning. To tackle this, mathematical modeling and optimization has become helpful tools for emergency responders and policy makers to evaluate required time for evacuation and to estimate numbers and distribution of casualties under different disaster scenario [11-14]. Such tools can help emergency professionals, and consequently the societies, to move toward building a more sustainable and resilient built environment [15,16].

One of the main challenges in evacuation planning is shelter allocation, i.e. the location of safe areas and the effective distribution of evacuees to these safe areas such that the total evacuation time is minimized. The process of selecting proper locations for safe areas can substantially affect the effectiveness of evacuation planning [17]. To address the uncertainties and the conflicting objectives, stochastic methods such as two-stage stochastic programming and stochastic scenario analysis are used to find the optimal shelter location and allocation. Sherali et al. [17] proposed a location-allocation model for a hurricane and flood scenario, and formulated the optimization problem as a nonlinear mixed integer programming problem. Balcik and Beamon [9] developed a linear programming model to determine the number and location of safe areas and the amount of relief supplies at each safe area. Rui et al. [18] formulated the transit evacuation operation as a location-routing problem (LRP) with uncertain demands, and solved the problem using hybrid genetic algorithms, artificial neural network, and hill climbing heuristic algorithms. He et al. [19] converted the work of Rui et al. into a multi-stage LRP, and resolved the problem using the same heuristic algorithms. Kongsomsaksakul et al. [20] formulated the optimal shelter location problem for a flood scenario as a Stackelberg game (bi-level programming), and solved it with Genetic Algorithm.

Although gradient based optimization techniques have specific applications in engineering optimization problems (see e.g. [21,22]), heuristic optimization methods have extensive applications when the search space is complex. To compare with gradient-based optimization techniques, heuristic methods can be more efficient in solving non-differentiable functions with many local optima [23]. Many heuristic methods have been developed in the past decades, such as Genetic Algorithm (GA) [24] which is based on the process of natural evolution, Ant Colony Optimization [25] which is based on the foraging behaviors of ant colonies, Particle Swarm Optimization [26] which is based on the interactions among a flock of birds, and Teaching- Learning-Based-Optimization (TLBO) [23] which is inspired by the learning process in classrooms where the influence of the teacher on the students and the interactions between students lead to a better overall performance of the class.

In this paper, the bi-level programming formulation for the shelter location-allocation problem is considered with a scenario-based approach to address the uncertainty in the population to be evacuated from a small town in Lombardy region, Italy. In this regard, four scenarios are developed to study the optimal number and location of shelters for evacuation during normal weekdays, at nights, during weekends, and during vacation times with visiting travelers. The rest of the paper is organized as follows: in section 2, the bi-level optimization formulation of the shelter location-allocation problem is presented; section 3 describes the solution methods to solve the upper- and lower-level optimization problems; in section 4, the study area is introduced along with the scenarios to be used; and in section 5, the results and their implications are presented.

## 2. Problem Formulation

The location-allocation of shelters for evacuation is formulated by Kongsomsaksakul et al. [20] as a bi-level programming problem to determine the location of shelters with capacity constraints. The problem is presented as a Stackelberg game where the number and location of shelters (or exit points) are determined by the authorities, as the leader, with the objective of minimizing the total network evacuation time, while evacuees, as the followers, decide how to minimize their travel time by selecting the shelter and the route to take given the location of the shelters. As a result of this formulation, the leader cannot directly control the behavior of the followers, but the leader's decision on the number and location of shelters will influence the behavior of followers.

### 2.1. Notation and Assumptions

The notation used in the bi-level location-allocation model is as follows:

$a$       index of links

$i$       index of origin nodes

Haghpanah and Foroughi. Optimal Shelter Location-Allocation| | |
|---|---|
| $j$ | index of potential shelter nodes |
| $\bar{j}$ | index of selected shelter nodes |
| $r$ | index of routes |
| $\beta$ | impedance parameter |
| $\delta_{ar}^{i\bar{j}}$ | link incidence parameter: 1, if link $a$ is on route $r$ between origin $i$ and shelter $\bar{j}$ ; 0, otherwise. |
| $f_r^{i\bar{j}}$ | traffic flow on route $r$ between origin $i$ and shelter $\bar{j}$ |
| $p_a$ | maximum acceptable degree of saturation for link $a$ |
| $q_{ij}$ | traffic demand between origin $i$ and shelter $j$ |
| $q_{i\bar{j}}$ | traffic demand between origin $i$ and shelter $\bar{j}$ |
| $V_a$ | traffic flow on link $a$ |
| $t_a$ | travel time on link $a$ |
| $C_a$ | traffic capacity of link $a$ |
| $K_j$ | capacity of shelter $j$ |
| $O_i$ | trip demand from zone $i$ |
| $R_{i\bar{j}}$ | route set for origin $i$ and shelter $\bar{j}$ |
| $X_j$ | 1, if shelter $j$ is selected; 0, otherwise. |

It is assumed that: (1) all the potential locations for shelters and their capacities are pre-determined; (2) the trip demand from each origin node is known; (3) only one type of vehicle (passenger car) is used in the network; and (4) the characteristics of the transportation network is given.

*2.2. Mathematical Formulation*

2.2.1. Upper Level

The upper level problem represents the optimization task of the authorities to minimize the total evacuation time on the network by selecting the number and location of the shelters under shelters' and links' capacity constraints. The upper level problem formulation is as follows:

$$\underset{X}{\text{Min}} \quad \sum_a V_a(X) t_a(V_a) \tag{1}$$

s.t.

$$\sum_i q_{ij} \leq K_j X_j \quad \forall j, \tag{2}$$

$$V_a(X) \leq p_a C_a \quad \forall a, \tag{3}$$

$$t_a = t_a^0 \left[ 1 + 0.15 \left( \frac{V_a}{C_a} \right)^4 \right], \tag{4}$$

$$X_j = \{0, 1\},$$

where Eq. (1) is the total evacuation time to be minimized, Eq. (2) is the shelter capacity constraint for each shelter $j$, Eq. (3) is the link capacity constraint for each link $a$, and Eq. (4) is the travel time function for link $a$, where $t_a^0$ is the free-flow travel time.

2.2.2. Lower Level

The lower level problem represents the optimization task of the evacuees to minimize self-evacuation time on the network by selecting the shelter and the route to exit the evacuation area. The lower level problem formulation is as follows:

$$\underset{V,q}{Min} \quad \sum_a \int_0^{V_a} t_a(w)dw + \frac{1}{\beta}\sum_i \sum_{\bar{j}} q_{i\bar{j}}\left(\ln q_{i\bar{j}} - 1\right) \tag{5}$$

s.t.

$$\sum_r f_r^{i\bar{j}} = q_{i\bar{j}} \qquad \forall i, \bar{j}, \tag{6}$$

$$\sum_{\bar{j}} q_{i\bar{j}} = O_i \qquad \forall i, \tag{7}$$

$$V_a = \sum_i \sum_{\bar{j}} \sum_r f_r^{i\bar{j}} \delta_{ar}^{i\bar{j}} \qquad \forall a, \tag{8}$$

$$t_a = t_a^0 \left[1 + 0.15\left(\frac{V_a}{C_a}\right)^4\right], \tag{9}$$

$$f_r^{i\bar{j}}, q_{i\bar{j}} \geq 0, \qquad \forall i, \bar{j}, r,$$

where Eq. (5) is the evacuees' evacuation time to be minimized, Eq. (6) is the flow conservation law for each origin-destination pair, Eq. (7) is the flow production law for each origin node, Eq. (8) defines the relation between routes, links, and nodes, and Eq. (9) is the travel time function for link $a$, where $t_a^0$ is the free-flow travel time.

3. Solution Method

Mainly, bi-level problems are not easy to solve since the solution of the upper level problem depends on the values of the lower level variables, which in turn, need the upper level variables to be determined. Furthermore, if the lower level problem includes non-linear constraints, the upper level problem can be considered as non-convex, which implies the ineffectiveness of standard optimization approaches for bi-level programming problems [27]. In this regard, as heuristic optimization methods can be used for any class of optimization problem, they have an extensive application in engineering problems, specifically bi-level optimization problems [28-30].

In this study, genetic algorithm is used to solve the bi-level shelter location-allocation problem described above. The upper level problem is solved using GA, and the lower level problem, i.e. the traffic assignment problem, is solved using the Double-Stage Algorithm [31] which is a simple algorithm for origin-destination trip estimation. To incorporate the constraints into the GA, the penalty method is used. The penalized objective function (OF) for the upper level problem is as follows:

$$OF_{upper-level} = \sum_a V_a(X) t_a(V_a) + \alpha \sum_j \max\left(q_{ij} - K_j X_j, 0\right) + \beta \sum_a \max\left(V_a - p_a C_a, 0\right), \tag{10}$$

where α and β are penalty parameters for constraints (2) and (3), respectively.

Figure 1 illustrates the overall scheme of the solution method. To address the stochasticity of the problem, the stochastic version of the Stackelberg-Nash-Cournot equilibrium model by De Wolf and Smeers [32] is used, which is basically based on scenario analysis.

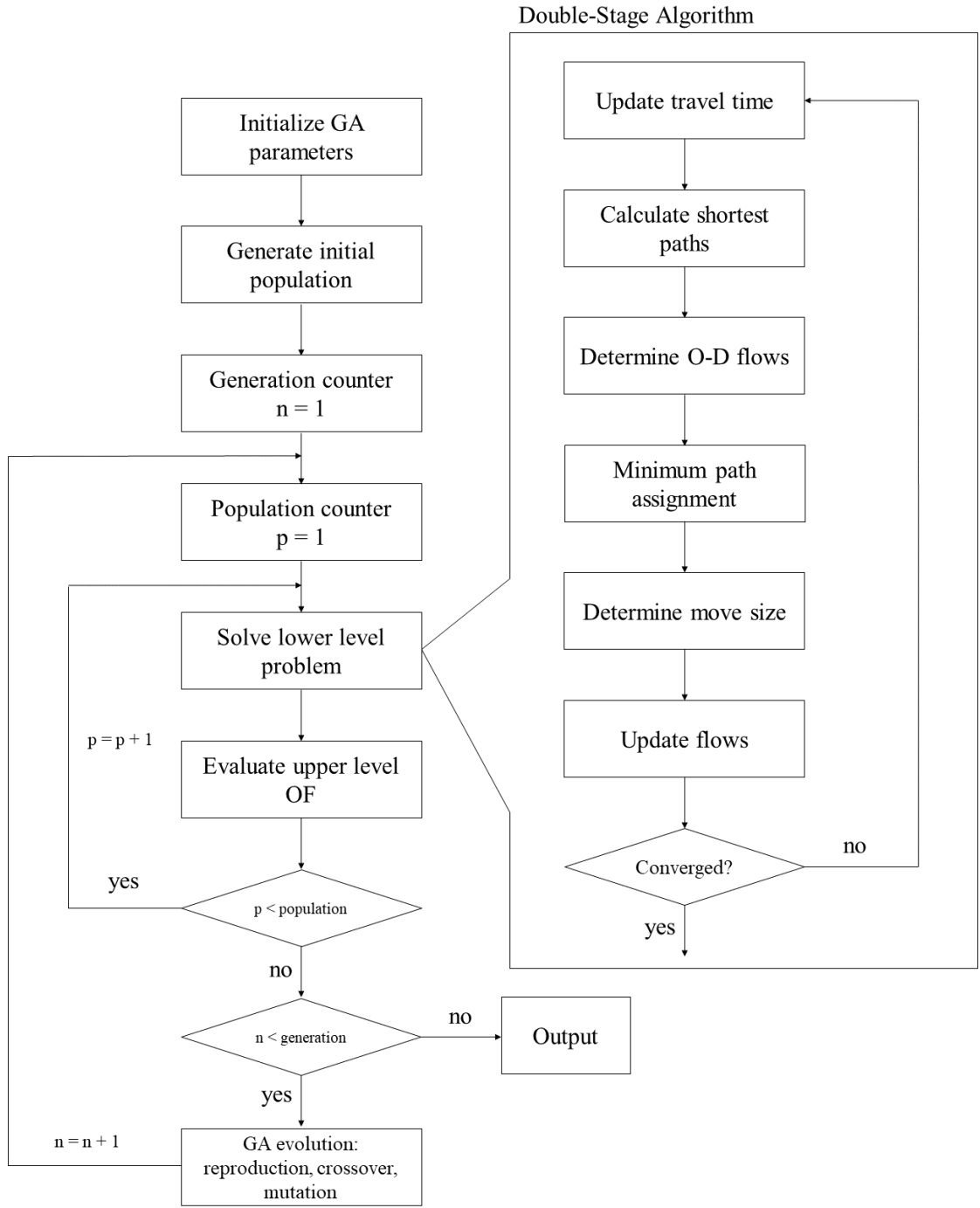

*Figure 1. Overall scheme of the solution method [20]*

**4. Study Area**

Po River passes through the boundary of Lombardy and Emilia-Romagna regions in Italy and is subjected to heavy flooding, specifically the southern part of the Lodi area adjacent to the left bank of the river which is located below the terrace morphology. This is a highly populated area that in case of flooding, would suffer substantial damage. The area is crossed by the most important arterials (roads and railways) linking the north and south of Italy. According to the flood risk maps, the township of San Rocco al Porto, located above the Po River, needs to be evacuated in case of a severe flood. San Rocco al Porto is a small town in the Province of Lodi in the Italian region Lombardy, located about 50 kilometers southeast of Milan and about 25 kilometers southeast of Lodi (see Figure 2).

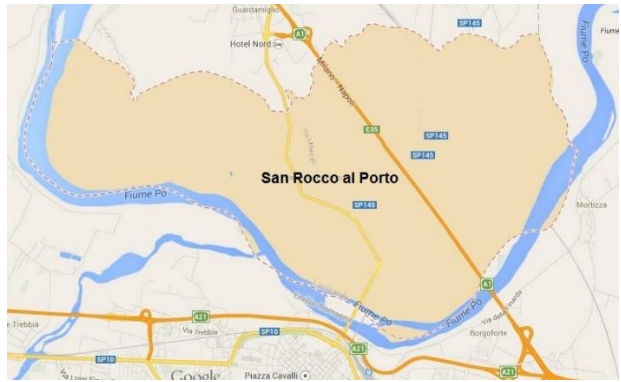

*Figure 2. Location of San Rocco al Porto with respect to the Po River (source: Google Map)*

*4.1. Demographic and Transportation Data*

A set of GIS shapefiles is provided including the administrative boundaries, municipalities and provinces, centroids, connectors, highway road network, local road network, flood risk regions, and location of variable message signs on highways. According to the ISTAT 2011 census, the San Rocco Al Porto town constitutes of 48 zones with a total population of 3250 (over 1300 households) from which 60% work within the town and thus, are considered for the evacuation. A total number of 2200 vehicles are registered in the city, from which 1300 vehicles will be used.

A total number of eight potential locations for shelters are determined as shown in Figure 3. Each potential shelter has a capacity of *1000 vph*. The capacity of main roads is set to *1000 vph* with a maximum acceptable saturation degree of 100% due to the emergency situation. The free-flow speed on the main roads is *35 mph*.

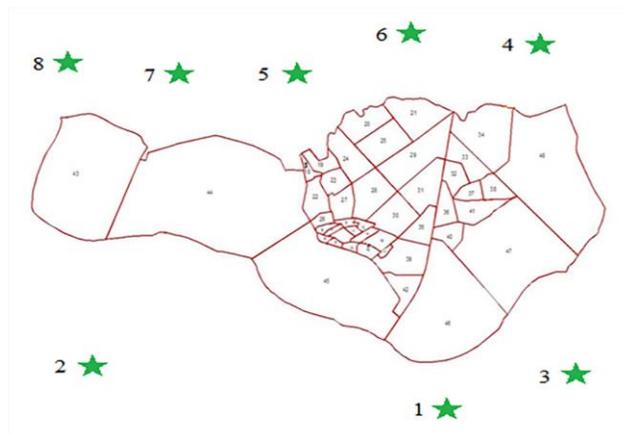

*Figure 3. Potential location of shelters*

*4.2. Scenarios*

The main sources of uncertainty in evacuation planning are the number and location of available shelters, trip production from origin nodes, and the characteristics of the transportation network. In this study, variable scenarios are set to study the optimal solution for different trip production rates from the origin nodes, which can be affected by the time of the event and demographic features of the region. Scenario 1 indicates the base scenario or the day scenario where most of the population is concentrated in the commercial areas, while a portion of the population is out of the town for work purposes. Scenario 2 is the night scenario where the majority of the population is located in the residential areas with 2200 vehicles. Scenario 3 is the weekend scenario which is similar to the night scenario except for a relatively more populated downtown, with a total of 2200 vehicles to be evacuated. Scenario 4 indicates the vacation time scenario where aside from the residents, there is a considerable number of visitors in the downtown area, increasing the number of vehicles to 4000.

## 5. Results and Discussion

The optimal solution for each scenario is obtained using the scheme depicted in Figure 1 with the following GA parameters: population size of 20, maximum number of generations of 50, 60% reproduction rate, and 40% probability of mutation. Moreover, the impedance parameter ($\beta$) is set to 10. In general, the higher the value of $\beta$, the more evacuees would be assigned to the nearest shelter. The results are shown in Table 1.

*Table 1. Optimal results of shelter location-allocation for the scenarios*

| Scenario | Shelters' attraction rates | | | | | | | | Min. travel time (vph) | Final evacuation time (min) |
|---|---|---|---|---|---|---|---|---|---|---|
| | 1 | 2 | 3 | 4 | 5 | 6 | 7 | 8 | | |
| 1 | 0 | 260 | 340 | 60 | 135 | 0 | 335 | 170 | 891.5 | 80 |
| 2 | 0 | 412 | 493 | 48 | 236 | 0 | 616 | 425 | 1494.0 | 80 |
| 3 | 0 | 390 | 575 | 90 | 203 | 0 | 605 | 340 | 1501.1 | 80 |
| 4 | 885 | 910 | 0 | 385 | 473 | 0 | 967 | 425 | 2984.4 | 85 |

According to the results, the day, night, and the weekend scenarios are similar where shelters 1 and 6 can be excluded, and total evacuation time will be 80 minutes. The total evacuation time is specifically dominated by the traffic flow to shelter 3. This is due to the fact that a considerable number of vehicles are located in areas which are only closed to shelter 3, and therefore, they all choose to exit the evacuation area only from this point. In scenario 4 (the vacation scenario), due to the extra number of visitors in the town, the traffic flows to shelters 1 to 3 (which are closer to the recreational areas) are increased. To prevent these shelters from oversaturation, the optimal solution to minimize the total travel time while meeting shelter and road capacity limits is to exclude shelter 3 so that the evacuees will select other shelters, and thus, the corresponding traffic will be distributed to other links and shelters.

Furthermore, the results of the scenarios can be used for more efficient resource allocation purposes. If the evacuation is to be ordered during day or night on non-vacation periods of the year, shelters 2, 3, and 7 will host the most evacuees; however, evacuation in periods with tourists visiting the town will bring more evacuees to shelters 1, 2, and 7.

**Conclusions**

The shelter allocation problem in evacuation, as one of the main challenges in evacuation planning, deals with finding the optimal location and number of shelters such that the total evacuation time would be minimized. There are different approaches to formulate such problem, each with advantages and shortcomings. The bi-level formulation of the problem has the advantage of addressing the fact that evacuees try to select shelters to which they can arrive the fastest; however, authorities would like to minimize the total evacuation time, i.e. the traffic flow on links are also important. Moreover, evacuees may not have information about the current capacity of shelters, or may not consider road capacity when selecting shelters, while shelter and road capacities are of the authorities' concern. These issues are reflected in the upper level and lower level objective functions and constraints. However, a shortcoming is neglecting information sharing during the evacuation, i.e. evacuees may decide to travel to another shelter due to oversaturation of some shelters.

Regarding the uncertainties in the shelter allocation problem, as the results show, variable parameters can lead to different solutions — a fact that planners and policy makers would like to have an estimation. Consequently, more stochasticity needs to be introduced to the problem, including road accessibility, evacuees' preparedness and response time, event duration, and event proximity. Certainly, obtaining the optimal shelter allocation settings for more scenarios with diverse characteristics can provide better tools for effective evacuation planning.